\documentclass[11pt]{amsart}   
\usepackage{euscript}   
\usepackage{amsfonts}
\usepackage{amsthm}
\usepackage{amsmath}
\usepackage{amsfonts}
\usepackage{latexsym}
\usepackage{amssymb}
\usepackage{mathtools}
\usepackage[usenames,dvipsnames]{color}
\usepackage[makeroom]{cancel}
\usepackage{tikz}
\usepackage{graphicx}

\begin{document}

\newtheorem{lemma}{Lemma}[section]
\newtheorem{theo}[lemma]{Theorem}
\newtheorem{coro}[lemma]{Corollary}
\newtheorem{rema}[lemma]{Remark}
\newtheorem{propos}[lemma]{Proposition}

\renewcommand{\thesection}{\arabic{section}}
\renewcommand{\theequation}{\thesection.\arabic{equation}}
\newcommand{\equnew}{\setcounter{equation}{0}}

\def\bdem{\begin{proof}}
\def\edem{\end{proof}}
\def\bequ{\begin{equation}}
\def\eequ{\end{equation}}

\newcommand{\hdos}{\mbox{$H^{2}$}}
\newcommand{\noi}{\noindent}
\newcommand{\ov}{\overline}
\newcommand{\papa}{\mbox{\(      H^{\infty}     \)}}
\newcommand{\rr}{\mbox{$    \longrightarrow   $}}
\newcommand{\clos}{\mbox{clos}}
\newcommand{\Om}{\Omega}
\renewcommand{\ker}{ \mbox{Ker}\, }
\newcommand{\la}{\langle}
\newcommand{\ra}{\rangle}
\newcommand{\disc}{ \mathbb{D} }
\newcommand{\N}{\mathbb{N}}
\newcommand{\re}{\mbox{Re}\,}
\newcommand{\bru}{\Delta}
\newcommand{\lamdos}{\lambda''_j}
\newcommand{\nuq}{q}

\newcommand{\matru}[4]{
\left[     \begin{array}{cc}
          #1     &     #2   \\
          #3     &     #4
       \end{array}
\right]                }

\hyphenation{cha-rac-te-ri-za-tion}  \hyphenation{res-pec-ti-ve-ly}
\hyphenation{ve-ri-fy}  \hyphenation{na-tu-ral}
\hyphenation{co-ro-lla-ry}  \hyphenation{en-cou-ra-ge-ment}
\hyphenation{o-ther-wise}  \hyphenation{re-gu-la-ri-zed}
\hyphenation{sa-tis-fies}  \hyphenation{pa-ra-me-ters}
\hyphenation{subs-tan-tia-lly}  \hyphenation{des-com-po-si-tion}
\hyphenation{sur-pri-sing}  \hyphenation{pro-ducts}  \hyphenation{cons-truc-tion}
\hyphenation{theo-rem}  \hyphenation{o-pe-ra-tor}  \hyphenation{o-pe-ra-tors} \hyphenation{se-pa-ra-ted}
\hyphenation{eigen-va-lues}  \hyphenation{dia-go-na-li-zable} \hyphenation{fo-llows}
\hyphenation{in-ter-po-la-ting}  \hyphenation{a-ppli-ed}  \hyphenation{mul-ti-pli-ci-ties}
\hyphenation{spe-ci-fying}


\title{Multi-orbital frames through model spaces}

\author[C. Cabrelli]{Carlos Cabrelli}
\address{The authors are at the Department of Mathematics, Facultad de Ciencias Exactas y Naturales,
Universidad de Buenos Aires.}
\author[U. Molter]{Ursula Molter}
\address{Carlos Cabrelli and Ursula Molter are also at: IMAS-CONICET-UBA, Argentina}
\email{carlos.cabrelli@gmail.com, umolter@conicet.gov.ar}
\author[D. Su\'arez]{Daniel Su\'arez}
\address{Daniel Su\'arez is also at IAM-CONICET, Argentina}
\email{dsuarez@dm.uba.ar}

\thanks{Supported by Universidad de Buenos Aires UBACyT 20020170100430BA, CONICET PIP11220150100355
and SECyT PICT 2014-1480}


\date{ \today }

\abstract
We characterize the normal operators $A$ on $\ell^2$ and the elements
$a^i \in \ell^2$, with $1\le i\le m$, such that the sequence
$$\{ A^n a^1 , \ldots , A^n a^m \}_{n\ge 0}$$
is a frame.
The characterization makes strong use of the pseudo-hyperbolic metric of $\disc$ and is given in terms
of the backward shift invariant subspaces of $H^2(\disc)$ associated to finite products
of interpolating Blaschke products.}

\maketitle

\section{Introduction}
\equnew

The study of the dynamical behaviour of a bounded operator $A$ on a Hilbert space $\mathcal{H}$ consists
of studying the orbits $\{ A^n f: \ n\in \mathbb{N}_0\}$ for $f\in \mathcal{H}$.
The literature is full of examples with characterizations of the operators $A$ such that
there exists an orbit satisfying a particular property, and sometimes also characterizations of the initial vector
for such orbits.  For instance, if $\mathcal{H}$ is separable and infinite dimensional,
the orbit $\{ A^n f \}$ is an orthonormal basis if and only if $A$ is the forward shift
with respect to the basis $e_n := A^nf$, for $n\ge 0$. Actually, this can be taken as the definition
of the forward shift with respect to a given ordered basis $\{ e_n\}_{n \ge 0}$. Moreover, it
is not hard to see that the only vectors with this property are $\lambda e_0$, where
$\lambda \in \mathbb{C}$, with $|\lambda|=1$.

Less restrictive requirements for an orbit is that of being a frame or even a Bessel sequence.
Motivated by a time-space sampling problem, in \cite{acmt} the authors characterize the
diagonalizable operators $A$ and the vectors $f\in \mathcal{H}$ such that the orbit
$\{ A^n f\}_{n\ge 0}$ is a frame for a Hilbert
space of numerable dimension. The problem is modeled with the space
$\mathcal{H} = \ell^2(\mathbb{N})$, but with the right definitions the result is valid for
finite dimension.

A fortiori, in \cite{accmp} it is shown that any normal operator $A$ that admits an orbit as a frame
must be diagonalizable, so the above result applies to normal operators as well.
The normality of $A$ allows the use of the spectral theorem, which
in conjunction with the fact that some orbit is a frame forces the operator to be diagonalizable.
Then, $\{ A^n f\}$ is a frame if and only if the sequence of eigenvalues $\{ \lambda_j \}$
is an interpolating sequence for the Hardy space of the disk $H^2(\disc)$, and
$$
f = \{ d_j (1-|\lambda_j|^2)^{\frac{1}{2}}: j\in \N \} \in \ell^2(\mathbb{N}),
$$
where $C^{-1} \le d_j \le C$ for all $j\in \mathbb{N}$ and some $C>0$
(see also \eqref{niunii} below and the subsequent comment).
The result holds for finite dimension, taking $\ell^2(J)$, where $J\subset \mathbb{N}$ is finite,
and accepting interpolating sequences also as those that perform finite interpolation in $H^2(\disc)$.

In \cite{cmpp} the authors consider the problem of characterizing the normal operators
$A$ and vectors $f_1, \ldots , f_m\in \ell^2(\mathbb{N})$, where $m\in \mathbb{N}$,
such that the union of orbits $\{ A^nf_j, \, n\in \mathbb{N}, \, 1\le j\le m \}$ is a frame.
They obtained a characterization where, as before, $A$ has to be diagonalizable, the eigenvalues form
a union of at most $m$ interpolating sequences for $H^2(\disc)$, and there are two more conditions,
the last of which is not well understood and difficult to handle.

In the present paper we give a different, more intuitive and geometric characterization,
which shows to what extent the pseudo-hyperbolic metric of $\disc$ plays a role in the structure
of the eigenvalues and their interaction with the vectors $f_j$.
To do so, we need some tools from the theory of $H^2(\disc)$, such as
interpolating sequences, reproducing kernels and model spaces, which we establish
in the next section.

Finally, in proving the above characterization we found a result of independent interest
(Thm.\@ \ref{difere}), which gives an upper bound for the Bessel constant of
the difference of normalized reproducing kernels in $H^2 = H^2(\disc)$
under some geometric conditions of their base points.

\subsection{The Hardy space $H^2$ and the model subspaces}

\noi
Write $\varphi_0(z)=\phi_0(z) =z$ and for $\lambda\neq 0$,
$$\varphi_{\lambda}(z) =   \frac{\lambda-z}{1-\ov{\lambda} z}
\, \ \ \mbox{ and }\ \ \,
\phi_{\lambda}(z) = \frac{\ov{\lambda}}{|\lambda_i|}  \,   \varphi_{\lambda}(z).
$$
If $\{\lambda_i\}$ is a sequence in $\disc$, the Blaschke product
$$ B(z) = \prod_{i} \phi_{\lambda_i}(z)  \mbox{ converges} \  \Leftrightarrow \ \,
\sum_i (1-|\lambda_i|^2) < \infty ,
$$
where the convergence is uniform on compact sets and $\{\lambda_i\}$ is called a Blaschke sequence.
Every function $f\in H^2$ factorizes as $f=gB$, where $g\in H^2$ has no zeros on $\disc$ and
$B$ is the Blaschke product of the zeros of $f$. If the zeros $\{ \lambda_j \}$ have single multiplicities,
the orthogonal complement $K_B :=(BH^2)^\bot$ of the
closed subspace $BH^2$ of $H^2$ is generated by
$$k_{\lambda_j}(z) = \frac{(1-|\lambda_j|^2)^{\frac{1}{2}}}{(1-\ov{\lambda}_jz)},
$$
the normalizations of the reproducing kernels $K_{\lambda_j}(z)=(1-\ov{\lambda}_jz)^{-1}$.
The name means that $\la f,K_{\lambda}\ra = f(\lambda)$ for every $f\in H^2$.
Blaschke products are special cases of inner functions, which are
functions $u\in H^\infty$ whose radial limit at the boundary satisfies $|u(e^{i\theta})|=1$
for almost every $e^{i\theta} \in \partial \disc$. The model spaces are $K_u := (uH^2)^\bot$,
which by Beurling's theorem \cite{beu} are the closed backward shift invariant subspaces of $H^2$.
They are called model spaces because the compression of the forward shift to $K_u$ is a model
for a broad class of contractions (see \cite{nfbk}).

A sequence $\{\lambda_j\}$ in $\disc$ is called interpolating (for $H^2$) if
$$E f := \{ \la f, k_{\lambda_j} \ra \} \in \ell^2, \ \, \forall f\in H^2
\ \mbox{ and every }\ w\in\ell^2 \ \mbox{ is of this form.}
$$
That is, $E: H^2 \to \ell^2$ is onto.
Here we allow the set of indexes of $\ell^2$ to be the set of natural numbers or a finite section, so finite sequences of different points
will also be called interpolating.

When the above holds, $\{ \lambda_i \}$ is the zero set of a Blaschke product $u$, and the restriction of
$E$ to the model space $E_u: K_u \to \ell^2$ is invertible. Therefore
$$
\| E_u^{-1} \|^{-2} \,\| f\|^2 \leq  \sum |\la f, k_{\lambda_i} \ra|^2 \leq  \| E_u \|^{2} \,\| f\|^2 ,
\ \hspace{0.4cm} \forall f\in K_u.
$$
This is equivalent to say that $\{ k_{\lambda_i} \}$ is a Riesz basis for $K_u$, or without specifying $u$,
that it is a Riesz sequence (see \cite[Lect.\@ 6, 1]{nik}).
This means that there are constants $C_0, \, C_1 >0$ such that
\bequ\label{dupy}
 C_0 \,\sum_j |c_j|^2 \leq  \|\sum_j c_j k_{\lambda_j} \|^2 \leq  C_1 \,\sum_j |c_j|^2 ,
\ \hspace{0.4cm} \forall \{ c_j\}\in \ell^2.
\eequ
On the other hand, $\{\lambda_j\}$ in $\disc$ is called interpolating (for $H^\infty$) if
$$\forall w\in \ell^\infty(J) \ \mbox{ there is }\ f\in H^\infty \ \mbox{ such that }\
f(\lambda_j) = w_j, \ \, \forall j\in J
$$
(again $J= \mathbb{N}$ or it is finite).
The problem of characterizing interpolating sequences for $H^\infty$ was considered
by several authors until Carleson obtained the definitive version in \cite{carin}.
In \cite{shsh} Shapiro and Shield provided a different proof and showed that
interpolating sequences are the same for all $H^p$, where $1\le p \le \infty$.
For a Blaschke sequence $\{ \lambda_i \}$ write $B$ for its Blaschke product and $B_j = \prod_{i: \, i\neq j} \phi_{\lambda_i}$.
The sequence is interpolating if and only if
$$
\delta(B) := \inf_j |B_j(\lambda_j)| >0 .
$$
A sequence satisfying this condition is usually called uniformly separated.
When this happens, $B_j(\lambda_j)^{-1} B_j\in H^2$,
with $\| B_j(\lambda_j)^{-1} B_j \| \le \delta^{-1}$ (here $\delta:= \delta(B)$), and
$$
f_i = k_{\lambda_i}  \ \mbox{ and } \   g_i = \frac{B_i}{B_i(\lambda_i)}  \,   k_{\lambda_i}
\ \mbox{ are biorthogonal sequences in $K_B$.}
$$
Therefore, when  $\{ c_i \} \in \ell^2$, the (unique) function $g\in K_B$ that interpolates
$\la g, k_{\lambda_j} \ra = c_j$ for all $j$ is $g(z) = \sum_i c_i g_i(z)$.
Any function $F\in H^2$  satisfying $\la F, k_{\lambda_j} \ra = c_j$ for all $j$ has the form
$F = g + Bh$, where $g$ is as above, $h\in H^2$, and $\|F\|^2 = \|g\|^2 + \|h\|^2$
(since $K_B= (BH^2)^\bot$ and multiplication by $B$ is an isometry).
In particular, $g\in K_B$ is the function of minimum norm that satisfies
$\la g, k_{\lambda_j} \ra = c_j$ for all $j$. Consequently,  \cite[Lemma 3]{shsh} gives us
\bequ\label{ddta}
\| g \|^2 \le (2/\delta^4) (1-2\log\delta) \sum |c_i|^2,
\eequ
where $\delta = \delta(B)$.
Also, the constants $C_0$ and $C_1$ of \eqref{dupy} depend only on $\delta$.
Indeed, a more general statement will be given for $C_1$ in Proposition \ref{carlbe}.
For $C_0$ notice that \eqref{ddta} together with Lemma  \ref{dualy} imply that
$$
\sum_i |\la f, g_i \ra|^2 \le C_\delta \|f \|^2, \hspace{7mm} \forall f\in K_{B},
$$
where $C_\delta$ is the constant of \eqref{ddta}.
In particular, when $f= \sum_j c_j k_{\lambda_j}$, for $\{ c_j \} \in \ell^2$,
we obtain the first inequality in \eqref{dupy} with $C_0= C^{-1}_\delta$.

The pseudo-hyperbolic metric in $\disc$ is given by $\rho(z,w)= |\varphi_z(w)|$, and
we denote the open ball
$$\bru(z,r) = \{ w\in \disc: \ \rho(z,w) <r \},\, \mbox{ where } 0<r<1,
$$
with the usual convention $\ov{\bru(z,r)}$  for the closed ball. Also, we will use that
Blaschke products satisfy the Lipschitz condition $\rho(B(z), B(w)) \le \rho(z,w)$ for $z,\, w\in \disc$,
and the elementary equality
\bequ\label{1menos}
1-|\varphi_v(z)|^2 = \frac{(1-|v|^2)(1-|z|^2)}{|1-\ov{v}z|^2}.
\eequ

\section{Basic necessary conditions}  
\noi
Let $\ell^2=\ell^2(J)$, where $J=\N$ or it is finite, and suppose that $A:\ell^2 \to \ell^2$ is a normal operator such that
there are $m$ vectors $a^1, \ldots , a^m\in\ell^2$ so that
$$
\mathcal{F}:= \{ A^n a^i:  \ n\in \N\cup\{0\}  , \  i=1, \ldots , m\}
$$
is a frame. If this happens, by exploiting the spectral theorem for normal operators it was shown in
\cite[Thm.\@ 5.6]{accmp} that $A$ is diagonalizable. So, from now on we assume that $A$ is a diagonal operator
with respect to the standard basis with eigenvalues $\{\lambda_j\}$. Next we aim to show some of the basic
properties that the $\lambda_j's$ and the vectors $a^i$ ($1\le i\le m$) must satisfy in order for
$\mathcal{F}$ to be a frame.

Let $e_{j_0}$ be the $j_0$ element of the standard basis and $a^i\in \ell^2$ for $i=1, \ldots, m$. Then
\begin{eqnarray*}
\sum_n \sum_{i=1} ^m |\la A^n a^i , e_{j_0} \ra|^2
&=& \sum_n  |\lambda_{j_0}^{2}|^{n} \Big[|a^1_{j_0}|^2 +\cdots + |a^m_{j_0}|^2\Big]  \\
&=& \frac{ |a^1_{j_0}|^2 +\cdots + |a^m_{j_0}|^2 }{1-|\lambda_{j_0}|^{2}  } .
\end{eqnarray*}
So, the lower bound for a frame implies that this expression is bounded below away from zero, implying that
$\sum_j (1-|\lambda_{j}|^{2}) \lesssim \sum_j \sum_{i=1}^m |a^i_{j}|^2 < \infty$, hence $\lambda_j$
is a Blaschke sequence.
Additionally, the Bessel constant (the upper frame constant) on the standard basis gives
$$
C_0 (1-|\lambda_{j}|^{2}) \le  \sum_{i=1} ^m |a^i_{j}|^2 \le  C_1 (1-|\lambda_{j}|^{2}) .
$$
In order to simplify notation it is convenient to consider a normalization $\tilde{a}^i$ of the vectors $a^i$.
For  $i=1,\ldots ,m\,$ write
$$\alpha^i_j = \ov{a^i_j} (|a^1_{j}|^2+\cdots +|a^m_{j}|^2)^{-\frac{1}{2}}
\ \ \mbox{ and }\ \
\tilde{a}^i_j = \ov{\alpha^i_j} (1-|\lambda_{j}|^{2})^{\frac{1}{2}}, \  (j\in J).
$$
Then $\sum_{i=1}^m|\alpha^i_j|^2 =1$ and  $a^i = d \cdot \tilde{a}^i$, a coordinate to coordinate product,
where $d\in \ell^\infty(J)$ is given by
$$
\sqrt{C}_0 \le d_j= \left[ \frac{|a^1_{j}|^2+\cdots +|a^m_{j}|^2}{1-|\lambda_{j}|^{2}}
\right]^{\frac{1}{2}} \le \sqrt{C}_1  .
$$
That is, any $m$ vectors $a^1,\ldots, a^m \in \ell^2(J)$ such that the union of the respective
$A$-orbits satisfies the lower and upper frame bounds when tested against the standard basis,
can be written as
\bequ\label{niunii}
a_j^i = d_j  \, \tilde{a}_j^i = d_j  \, \ov{\alpha^i_j} (1-|\lambda_{j}|^{2})^{\frac{1}{2}},
\ \ \mbox{ for $j\in J\,$ and $\,1\le i\le m$},
\eequ
where $C^{-1} \le d_j\le C\,$  for some $\,C\ge 1\,$ and $\,\sum_{i=1}^m|\alpha^i_j|^2 =1\,$
for all $j\in J$.
Moreover, it is clear that
$\{ A^n a^i \!: \, n\ge 0, \, 1\le i\le m \}$ is a frame (a Bessel sequence)
if and only if  $a^i$ are  given by \eqref{niunii} and
$\{ A^n \tilde{a}^i \!: \, n\ge 0, \, 1\le i\le m \}$ is a frame (respectively, a Bessel sequence).
So, from now on we work with $\tilde{a}^i$ for $1\le i\le m$.\vspace{1mm}

Write $\tilde{b} = \{ (1-|\lambda_{j}|^{2})^{\frac{1}{2}} \} \in \ell^2(J)$,
and let $c\in \ell^2(J)$, where $J$ could be finite. Then
\begin{eqnarray*}
\sum_n |\la A^n \tilde{b} , c \ra|^2 &=&
\sum_n \sum_{i,j} \lambda_i^n \ov{\lambda}_j^n (1-|\lambda_i|^2)^{\frac{1}{2}} (1-|\lambda_j|^2)^{\frac{1}{2}} \,
\ov{c}_i c_j \nonumber \\
 &=&
\sum_{i,j} \sum_n \lambda_i^n \ov{\lambda}_j^n (1-|\lambda_i|^2)^{\frac{1}{2}} (1-|\lambda_j|^2)^{\frac{1}{2}} \,
\ov{c}_i c_j \ \nonumber \\
 &=&
\sum_{i,j} \frac{ (1-|\lambda_i|^2)^{\frac{1}{2}} (1-|\lambda_j|^2)^{\frac{1}{2}} }{1- \lambda_i \ov{\lambda}_j} \,
\ov{c}_i c_j
\nonumber \\
 &=&  \sum_{i} \sum_{j} \la  c_j  k_{\lambda_j}, c_i k_{\lambda_i}\ra
\ =  \  \|   \sum_{j}c_j k_{\lambda_j}\|^2.
\end{eqnarray*}
It follows that
\bequ\label{kcat}
\ \ \ \ \sum_n |\la A^n \tilde{a}^i , c \ra|^2 =  \|   \sum_{j}\alpha^i_j c_j k_{\lambda_j}\|^2
\ \ \ \mbox{ for $1\le i\le m$}.
\eequ
In particular, when $m=1$, the orbit $\{ A^n \tilde{a}^1, n\ge0\}$ is a frame for $\ell^2$ if and only if
$\{k_{\lambda_j} \}$ is
a Riesz basis for the subspace $K_u =(uH^2)^\bot$, where $u$ is the Blaschke product with zeros $\lambda_j$.
By the previous section this happens if and only if $\{ \lambda_j\}$ is an interpolating sequence
(see  \cite[Thm.\@ 3.14]{acmt}).

\subsection{Carleson measures and Bessel sequences}

A positive measure $\mu$ on $\disc$ is called a Carleson measure if
$$\int |f|^2 d\mu \le C^2_2 \|f\|^2\ \ \ \ \forall \  f\in H^2.
$$
It is well known  (see  \cite[I, Thm.\@5.6]{gar}) that $\mu$ is Carleson if and only if
$$\mu(Q) \le C \ell(Q)$$
for every angular square $Q= \{ re^{i\theta}:  1-\ell \le r <1, \ |\theta -\theta_0| \le \ell \} ,$
where $\ell = \ell(Q)$.
The smallest constant $C$ is called the Carleson norm of $\mu$ and  is
denoted by $\| \mu \|_\ast$. Also, the optimal $C_2$
and  $\| \mu \|_\ast$ are equivalent quantities.

Carleson measures and interpolating sequences are closely related, as we shall see in the next lemma.
For a sequence $\{ \lambda_j \}$ in $\disc$ consider the purely atomic measure
$\mu= \sum_j (1-|\lambda_j|^2) \delta_{\lambda_j}$, where
$\delta_\lambda$ is the probability measure with mass concentrated at $\lambda\in \disc$.
Observe that $\mu(\disc) = \sum_j (1-|\lambda_j|^2) <\infty$ if and only if
$\{\lambda_j\}$ is a Blaschke sequence. Furthermore,  it is well known that the following holds.
\begin{lemma}\label{becar}
Let\/ $S= \{ \lambda_j \}$ be a sequence in $\disc$ and $\mu= \sum_j (1-|\lambda_j|^2) \delta_{\lambda_j}$.
Then
\begin{enumerate}
\item[(1)] $S$ is a finite union of interpolating sequences if and only if $\mu$ is Carleson.
\item[(2)] $S$ is interpolating if and only if $\mu$ is Carleson and
$\rho (\lambda_j , \lambda_k) \ge \beta >0$ when $j\neq k$ (i.e.: $S$ is separated).
\end{enumerate}
When {\em (2)} holds and $B$ is the respective Blaschke product, $\delta(B)$ can
be estimated from $\mu$ and $\beta$, and vice versa.
\end{lemma}
\bdem
Assertion (1) is proved in \cite[Lemma 21]{mc-s}. Assertion (2) can be found in \cite[VII, Thm.\@1.1]{gar}.
The same theorem shows the equivalence of (2) with $S$ being uniformly separated and the
relations between the various parameters are established.
\edem

\noi The following basic and well-known result on Bessel sequences can be found, for instance,
in  \cite[pp.\@ 51-53]{ole}.

\begin{lemma}\label{dualy}
For a sequence $\{ f_j \}$  in $\mathcal{H}$, the next assertions are equivalent:

\begin{enumerate}
\item[(1)] $T (\{ c_j \}) = \sum c_j f_j$ is a bounded operator from $\ell^2$ to $\mathcal{H}$:
$$\| \sum c_j f_j\|^2 \le \| T \|^2 \, \sum |c_j|^2.
$$
\item[(2)] $T^\ast f =  \{ \la f, f_j\ra  \}$ is a bounded operator from $\mathcal{H}$ to $\ell^2$:
$$
\sum |\la f, f_j\ra|^2 \le \| T^\ast \|^2 \, \|f\|^2.
$$
\end{enumerate}
\end{lemma}

\noi
In \cite{phi} Philipp uses the notion of Carleson measure on spectral measures to characterize
orbits of normal operators that are Bessel sequences. This connection becomes particularly
clear when dealing with reproducing kernels, as our next result shows.
First, notice that
by \eqref{kcat}, $\{ A^n \tilde{a}^1 , \ldots , A^n \tilde{a}^m \}$
is a Bessel sequence with constant $\le B^2$ if and only if
$\{ \alpha^1_j k_{\lambda_j}  ,\ldots ,\alpha^m_j k_{\lambda_j} \}$  is Bessel with constant $\le B^2$.

\begin{propos}\label{carlbe}
The sequence $\{ A^n \tilde{a}^1 , \ldots , A^n \tilde{a}^m \}$ is Bessel with constant $\le B^2$
if and only if the measure  $\mu:= \sum_{i}  (1-|\lambda_i|^2) \delta_{\lambda_i}$ is Carleson,
where if\/ $C\ge 0$ is such that
$$\int |f|^2 d\mu \le C^2 \|f\|^2, \ \ \ \mbox{ for }\ f\in H^2,$$
then we can take $B^2 \le C^2 \le m B^2$.
\end{propos}
\bdem
If
\bequ\label{nuy}
 \big\| \sum_{j}\alpha^1_j c_j k_{\lambda_j}\big\|^2  +\cdots +
 \big\| \sum_{j}\alpha^m_j c_j k_{\lambda_j}\big\|^2   \le   B^2 \, \sum_{j}|c_j|^2,
\eequ
the inequality holds for each member of the left sum.
Hence,  Lemma \ref{dualy} says that for all $f\in H^2$,
\begin{align*}
\sum_{j}\sum_{i=1}^m |\alpha^i_j|^2 \, |\la f,   k_{\lambda_j}\ra |^2
&= \sum_{i=1}^m \sum_{j}  |\la f, \alpha^i_j  k_{\lambda_j}\ra |^2
\ \le \   B^2 \sum_{i=1}^m\|f\|^2.
\end{align*}
Since $|\alpha^1_j|^2+\cdots+|\alpha^m_j|^2=1$, then
$\mu:= \sum_{j}  (1-|\lambda_j|^2) \delta_{\lambda_j}$ is Carleson with
$$\int |f|^2 d\mu \le mB^2 \|f\|^2,  \ \ \ \ \mbox{ for $f\in H^2$.}
$$

\noi
Reciprocally, if $\mu$ is Carleson with $\sum_{j}|\la f,  k_{\lambda_j}\ra |^2 \le  C^2  \|f\|^2$,
by Lemma \ref{dualy},
$$\big\| \sum_{j}\alpha^i_j c_j k_{\lambda_j}\big\|^2 \le C^2 \, \sum_{j}|\alpha^i_j|^2  |c_j|^2
 \ \ \ \mbox{ for $1\le i\le m\,$}
$$
and $\{ c_j \} \in \ell^2$. Adding for $1\le i\le m$ we get \eqref{nuy} with $B^2 = C^2$.
\edem

\noi The following estimate for the distance between two normalized reproducing kernels in $H^2$
is sharp and can be found in \cite[Lemma\@ B7]{cmpp}:
\bequ\label{labubu}
\| k_v-k_w \|^2 \le 2 \rho(v,w)^2 \hspace{7mm} \forall v,w \in \disc .
\eequ

\vspace{0.1mm} 
\begin{theo}\label{kasem}
Let $S = \{\lambda_j \}$ be a sequence such that $\sum_{i}  (1-|\lambda_i|^2) \delta_{\lambda_i}$ is Carleson
and let $m\in \mathbb{N}$.
Then there is a constant $D>0$ satisfying
\bequ\label{mini}
 D^2 \, \sum_{j}   |c_j|^2 \le
 \sum_{i=1}^m   \big\|   \sum_{j}\alpha^i_j c_j k_{\lambda_j}\big\|^2
 \hspace{6mm} \ \forall c\in \ell^2 .\vspace{-2mm}
\eequ
if and only if  there is\/ $\eta>0$ such that
\begin{enumerate}
\item[(1)]  $\bru(\lambda_j , \eta)$ contains no more than m points of $S$
counting repetitions for all $j$.
\item[(2)]
if\/ $\lambda_{j_1} , \ldots, \lambda_{j_p}$ ($p\le m$) are the points of $S$ in
$\bru(\lambda_{j_1}, \eta)$  counting repetitions, the related matrix satisfies
$$
D_0^2
\left\|
    \begin{bmatrix}
    c_{j_1} \\
    \vdots  \\
    c_{j_p} \\
  \end{bmatrix}
\right\|^2_{\mathbb{C}^p}
\le
\left\|
  \begin{bmatrix}
    \alpha^1_{j_1}  & ...  &\alpha^1_{j_p} \\
    \vdots     &            &  \vdots  \\
     \alpha^m_{j_1}  & ...  &\alpha^m_{j_p}\\
  \end{bmatrix}
  \begin{bmatrix}
    c_{j_1} \\
    \vdots \\
    c_{j_p} \\
  \end{bmatrix}
\right\|^2_{\mathbb{C}^{m \times 1}}
\ \ \ \forall (c_{j_1},\ldots, c_{j_p}) \in \mathbb{C}^p,
$$
where $D_0>0$ does not depend on $p$ or the $\alpha$'s .
\end{enumerate}
\end{theo}

\bdem[Proof of necessity for Theorem \ref{kasem}]
First we show that there exist an $\eta$ such that (1) holds.
Suppose otherwise that for any $\eta >0$ there are at least m+1 points
$\lambda_{j_0}, \lambda_{j_1},\ldots ,  \lambda_{j_m}$ of $S$ counting repetitions, such that
\bequ\label{etu1}
\lambda_{j_0}, \lambda_{j_1}, \ldots , \lambda_{j_m}\in \bru(\lambda_{j_0}, \eta).
\eequ
To simplify notation we assume that $j_s = s$ for $s=0,\ldots,m$.
Taking  $c= (c_0, \ldots, c_m , 0 , \ldots ) \in \ell^2$, each summand of
the right side of \eqref{mini} is
\begin{eqnarray*}
\lefteqn{ \hspace{3mm} \big\| \sum_{j=0}^m\alpha^i_j c_j k_{\lambda_j}\big\|^2 = }\\*[-1mm]   &=&
 \Big\| \Big[\sum_{j=0}^m   c_j \alpha^i_j\Big] k_{\lambda_0}  + c_1 \alpha^i_1 [k_{\lambda_1}-k_{\lambda_0} ]
+ \cdots + c_m \alpha^i_m [k_{\lambda_m}-k_{\lambda_0} ]\Big\|^2 ,
\end{eqnarray*}
where $1\le i\le m$. If we take a normalized vector $c\in \mathbb{C}^{m+1}$ such that
\begin{equation*}
  \begin{bmatrix}
    \alpha^1_0 &      \ldots & \alpha^1_m \\
     \vdots    &             & \vdots   \\
    \alpha^{m}_0 &  \ldots & \alpha^{m}_m
  \end{bmatrix}
  \begin{bmatrix}
    c_0 \\
    \vdots \\
    c_m \\
  \end{bmatrix}
   =
  \begin{bmatrix}
    0 \\
    \vdots \\
    0 \\
  \end{bmatrix},
\end{equation*}
by the Cauchy-Schwarz inequality and \eqref{labubu},  the right side of \eqref{mini} becomes
\begin{eqnarray*}
\lefteqn{  \sum_{i=1}^m   \big\|  c_1 \alpha^i_1 [k_{\lambda_1}-k_{\lambda_0} ]
+ \cdots + c_m \alpha^i_m [k_{\lambda_m}-k_{\lambda_0} ]  \big\|^2  \le } \hspace{15mm}  \\*[-2mm]
&\le&
 m \sum_{i=1}^m \Big( \big\|  c_1 \alpha^i_1 [k_{\lambda_1}-k_{\lambda_0} ]\big\|^2
+ \cdots + \big\|c_m \alpha^i_m [k_{\lambda_m}-k_{\lambda_0} ]  \big\|^2  \Big)    \\
&\le&
 m 2\eta^2  \sum_{i=1}^m  \big(|c_1 \alpha^i_1|^2 + \cdots + |c_m \alpha^i_m|^2 \big) \le 2m\eta^2.
\end{eqnarray*}
Therefore, \eqref{mini} applied to this particular case says that $D^2  \le 2m \eta^2$.
This means that  \eqref{etu1} can't happen for $\eta < D/\sqrt{2m}$.\\

\noi
Assume now that  $\eta$ satisfies (1)
and suppose that
for $1\le p\le m$, we have
$$
\lambda_{j_1}, \ldots , \lambda_{j_p}\in \bru(\lambda_{j_1}, \eta).
$$
As before, we write $j=1, \ldots , p$ instead of $j_1, \ldots, j_p$. Then
\begin{eqnarray*}
 \big\| \sum_{j=1}^p\alpha^i_j c_j k_{\lambda_j}\big\|^2   \!  &=&  \!
 \Big\| \Big[\sum_{j=1}^p   c_j \alpha^i_j\Big] k_{\lambda_1}  + c_1 \alpha^i_1 [k_{\lambda_2}-k_{\lambda_1} ]
+ \cdots + c_p \alpha^i_p [k_{\lambda_p}-k_{\lambda_1} ] \Big\|^2 \\
&\hspace{-47mm}\le& \hspace{-25mm}   m\left[
  \Big\| \Big[\sum_{j=1}^p   c_j \alpha^i_j\Big] k_{\lambda_1}  \Big\|^2 +
 \| c_1 \alpha^i_1 [k_{\lambda_2}-k_{\lambda_1} ] \big\|^2
+  \cdots + \| c_p \alpha^i_p [k_{\lambda_p}-k_{\lambda_1} ]\|^2  \right] .
\end{eqnarray*}

\noi So, if  $c= (c_1, \ldots, c_p , 0 , \ldots ) \in \ell^2$, by \eqref{mini} applied to this case
and \eqref{labubu},
\begin{eqnarray*}
 D^2 \sum_{j=1}^p  |c_j|^2     &\le &    \sum_{i=0}^m
 \big\| \sum_{j=1}^p\alpha^i_j c_j k_{\lambda_j}\big\|^2  \\
&\le&
\sum_{i=0}^m       m\left(
  \Big|\sum_{j=1}^p   c_j \alpha^i_j\Big|^2 + 2 \eta^2  |c_1 \alpha^i_1|^2
+  \cdots + 2 \eta^2  |c_p \alpha^i_p|^2 \right)  ,
\end{eqnarray*}
which clearly implies that
$$
D^2
\left\|
    \begin{bmatrix}
    c_{ 1} \\
    \vdots  \\
    c_{ p} \\
  \end{bmatrix}
\right\|^2_{\mathbb{C}^p}
\le
m\left\|
  \begin{bmatrix}
    \alpha^1_{ 1}  & ...  &\alpha^1_{ p} \\
    \vdots     &            &  \vdots  \\
     \alpha^m_{ 1}  & ...  &\alpha^m_{ p}\\
  \end{bmatrix}
  \begin{bmatrix}
    c_{ 1} \\
    \vdots \\
    c_{ p} \\
  \end{bmatrix}
\right\|^2_{\mathbb{C}^{m \times 1}}
\! \! \! +2m\eta^2
\left\|
    \begin{bmatrix}
    c_{ 1} \\
    \vdots  \\
    c_{ p} \\
  \end{bmatrix}
\right\|^2_{\mathbb{C}^p}
$$
for every $(c_{ 1},\ldots, c_{ p}) \in \mathbb{C}^p$.
Thus, if $\Upsilon$ denotes the above matrix and $c\in \mathbb{C}^p$ is normalized,
when $2m \eta^2 < D^2$ we get
$$
0<D_0^2=\frac{D^2- 2m \eta^2}{m} \le \| \Upsilon c\|^2 .
$$
\edem


\subsection{Differences of normalized reproducing kernels}

For $E\subset \disc$ and $0<r<1$ we write $\Om_r(E):= \{ z\in\disc: \rho(z,E)\le r \}$.
\begin{lemma}    
Let $\mu_0= \sum (1-|\lambda_j|^2) \delta_{\lambda_j}$ be a Carleson measure, $0<r<1$
and\/ $\lambda'_j \in \bru(\lambda_j, r)$.
Then $\mu_r= \sum (1-|\lambda'_j|^2) \delta_{\lambda'_j}$ is a Carleson measure
such that for some constant $C(r) \ge 1$,
\bequ\label{omy}
C(r)^{-1} \,  \|\mu_r\|_{\ast}  \le  \|\mu_0\|_{\ast} \le C(r) \,  \|\mu_r\|_{\ast} .
\eequ
\end{lemma}
\bdem
By \cite[p.$\,$3]{gar} any $z\in \ov{\bru(\lambda_j, r)}$ satisfies $|z|  \le  \varphi_{|\lambda_j|} (-r)$.
Then \eqref{1menos} implies that for any angular square $Q\subset \disc$,
$$
\mu_r(\Om_r(Q))  \ge  \sum_{\lambda_j \in Q} (1-|\varphi_{|\lambda_j|} (-r)|^2)
=  \sum_{\lambda_j \in Q}  \frac{(1-|\lambda_j|^2)(1-r^2)}{(1 + |\lambda_j| r)^2 }
\ge  \frac{1-r}{1+r}   \mu_0(Q) .
$$
If $Q_r$ is the smallest angular square containing $\Om_r(Q)$, there is $c(r) \ge 1$ such that
$\ell(Q_r) \le c(r) \ell(Q)$. Hence,
$$
 \left[\frac{1-r}{1+r} \right] \,   \frac{\mu_0(Q)}{c(r)\ell(Q)} \le  \frac{\mu_r(Q_r)}{\ell(Q_r)} ,
$$
implying that
$$
  \|\mu_0\|_{\ast} \le \frac{1+r}{1-r} c(r) \,  \|\mu_r\|_{\ast} = C(r) \,  \|\mu_r\|_{\ast}.
$$
Hence the lemma follows by symmetry.
\edem

\vspace{1mm}
\noi In what follows  a Bessel sequence $\{ f_j \}$ in $H^2$, that is,
$$
\sum_j |\la f ,f_j\ra |^2 \le B^2 \|f\|^2   \ \ \ \ \forall f\in H^2, \vspace{-3mm}
$$
we write $\mathcal{B}^2(\{ f_j\})$ for the smallest constant $B^2$.\\

\begin{coro}\label{muss}
Let $\mu_0= \sum (1-|\lambda_j|^2) \delta_{\lambda_j}$ be a Carleson measure,  $0<r<1$
and $\lambda'_j \in \ov{\bru(\lambda_j, r)}$.
Then $\{ k_{\lambda'_j} \}$ is Bessel with
$\mathcal{B}^2(  \{ k_{\lambda'_j}\} ) \le C(r, \|\mu_0\|_\ast)$.
\end{coro}
\bdem
By \eqref{omy} $\mu_r= \sum (1-|\lambda'_j|^2) \delta_{\lambda'_j}$   is a Carleson
measure with $\| \mu_r \|_\ast$ depending on $r$ and $\| \mu_0 \|_\ast$.
Therefore,  the comments preceding Lemma \ref{becar} say that
there is  $C_r\ge 0$ depending on $\| \mu_r \|_\ast$  such that
$$\int |f|^2 d\mu_r \le C_r^2 \|f\|^2, \ \ \ \mbox{ for }\ f\in H^2.$$
Thus,  $\{ k_{\lambda'_j} \}$ is Bessel with
$\mathcal{B}\big(\{ k_{\lambda'_j} \}\big)\le C_r^2$ by Prop.$\ $\ref{carlbe}.
\edem

\vspace{1.5mm}
\noi
If $\lambda_1, \ldots , \lambda_N\in \disc$ and  $\lambda'_j \in \bru( \lambda_j, \eta)$
for some $0<\eta <1$, then \eqref{labubu} implies that the Bessel constant
$$
\mathcal{B}^2\big(  \{  k_{\lambda_j}-k_{\lambda'_j}  \}\big) \le 2N  \eta^2.
$$
Next we see how to control this constant for infinitely many values of $\lambda_j$.
Together with Theorem \ref{kasem}, this is the main result of the paper.

\begin{theo}\label{difere}
Let $\mu_0= \sum (1-|\lambda_j|^2) \delta_{\lambda_j}$ be a Carleson measure, $0<r<1$,
and  $\lambda'_j \in \bru( \lambda_j, \eta)$, where\/ $0<\eta <r$.
Then there is a constant $C>0$ depending only on $r$ and $\{ \lambda_j \}$
such that
$$
\mathcal{B}\big(\{   k_{\lambda_i}-k_{\lambda'_i}\}  \big) \le C  \eta.
$$
\end{theo}
\bdem
Since by Lemma \ref{becar}, $\{\lambda_i \}$ is a finite union of interpolating sequences,
we can assume that it is interpolating.
Let $B$ be the Blaschke product with zeros $\lambda_i$, write $B_i= B/ \phi_{\lambda_i}$
(i.e.: $B$ with the factor $\phi_{\lambda_i}$ removed) and recall that
$\delta(B) = \inf_i |B_i(\lambda_i)|>0$.
We prove first the result for $f\in K_B$, which by \eqref{ddta} can be written as
\bequ\label{tupi}
f= \sum_i c_i \frac{B_i}{B_i(\lambda_i)} k_{\lambda_i},
\ \mbox{ with $c\in \ell^2\ $ and }\  \|f\| \le C_\delta \|c\|_{\ell^2} ,
\eequ
where $C_\delta>0$ is a constant depending only on $\delta(B)$.
So,
\begin{align}
\la f , k_{\lambda_j}-k_{\lambda'_j} \ra &=
\Big\la c_j \frac{B_j}{B_j(\lambda_j)} k_{\lambda_j} , k_{\lambda_j}-k_{\lambda'_j} \Big\ra +
\Big\la \sum_{i: i\ne j} c_i \frac{B_i}{B_i(\lambda_i)} k_{\lambda_i}, k_{\lambda_j}-k_{\lambda'_j} \Big\ra \nonumber\\
&= c_j \left(1- \frac{B_j(\lambda'_j)}{B_j(\lambda_j)} \la k_{\lambda_j} , k_{\lambda'_j} \ra \right)-
\Big\la \sum_{i: i\ne j} c_i \frac{B_i}{B_i(\lambda_i)} k_{\lambda_i}, k_{\lambda'_j} \Big\ra \nonumber\\
&=  D_j -R_j . \ \label{qu1}
\end{align}
$$
\hspace{-12mm}\mbox{To estimate } \ \,
D_j =
\frac{c_j}{B_j(\lambda_j)}
\left[B_j(\lambda_j)-  B_j(\lambda'_j) +
B_j(\lambda'_j) ( 1- \la k_{\lambda_j} , k_{\lambda'_j} \ra ) \right]
$$
we notice that
\bequ\label{qu2}
|D_j| \le
\left|\frac{c_j}{B_j(\lambda_j)}\right| \,
\left[   2\rho(\lambda_j , \lambda'_j) + \sqrt{2}\rho(\lambda_j , \lambda'_j) \right]
\le  |c_j| \, 4 \frac{\rho(\lambda_j , \lambda'_j)}{\delta(B)},
\eequ
where the first inequality comes from
$\rho(B_j(\lambda_j),  B_j(\lambda'_j))  \le \rho(\lambda_j,  \lambda'_j)\,$
and from \eqref{labubu}.

To estimate $R_j$ write $B_{i,j}= B/ (\phi_{\lambda_i}\phi_{\lambda_j})$.
For $0<r<1$ consider the analytic function $F_j: \ov{\bru(\lambda_j , r)} \to \mathbb{C}$ given by
$$
F_j(\lambda') =
 \sum_{i: i\ne j} \frac{c_i}{B_i(\lambda_i)}\,  \la B_{i,j} k_{\lambda_i},  K_{\lambda'} \ra ,
$$
where $K_\lambda(z) = (1-\ov{\lambda}z)^{-1}$ is the reproducing kernel for $H^2$.
By the maximum modulus principle $F_j$ attains its maximum on the boundary of   $ \bru(\lambda_j , r)$.
That is, there is $\lamdos \in \partial\bru(\lambda_j , r)$ (i.e.: $\rho(\lambda_j , \lamdos)=r$) such that
$$
|F_j(\lambda')| \le |F_j(\lamdos)| \ \mbox{ for every }\  \lambda' \in \ov{\bru(\lambda_j , r)}.
$$
Since $\lambda' = \varphi_{\lambda_j}(w)$ with $0\le |w| \le r$,
a straightforward estimate from formula \eqref{1menos} gives
$$
\left[\frac{1-r}{1+r}\right]   (1-|\lambda_j|^2)
\le (1-|\lambda'|^2)  \le
\left[\frac{1+r}{1-r}\right]    (1-|\lambda_j|^2) ,
$$
implying that
\bequ\label{atla}
(1-|\lambda'|^2) |F_j(\lambda')|^2 \le C_1(r)\, (1-|\lamdos|^2) |F_j(\lamdos)|^2
\eequ
for some constant $C_1(r) >0$. Since
$|\phi_{\lambda_j} (\lamdos)| = \rho(\lambda_j, \lamdos) =r$, then
\begin{align*}
r(1-|\lamdos|^2)^{\frac{1}{2}} |F_j(\lamdos)| &=
\Big| \phi_{\lambda_j} (\lamdos)\sum_{i: i\ne j} \frac{c_i}{B_i(\lambda_i)}\, \la B_{i,j} k_{\lambda_i},  k_{\lamdos} \ra \Big| \\
&=  \Big| \sum_{i: i\ne j} \frac{c_i}{B_i(\lambda_i)}\,
\la B_{i} k_{\lambda_i},  k_{\lamdos} \ra \Big| \\
&\le  \Big| \sum_{i} \frac{c_i}{B_i(\lambda_i)}\,
\la B_{i} k_{\lambda_i},  k_{\lamdos} \ra \Big| +   \Big|\frac{c_j}{B_j(\lambda_j)}\,
\la B_{j} k_{\lambda_j},  k_{\lamdos}\ra \Big|.
\end{align*}
Consequently, the Cauchy-Schwarz inequality gives
\begin{align}
r^{2}\sum_j (1-|\lamdos|^2) |F_j(\lamdos)|^{2}
&\le 2\sum_j |\la f, k_{\lamdos}\ra|^2 + 2\sum_j \Big|\frac{c_j}{B_j(\lambda_j)}\Big|^2 \nonumber\\
&\le 2\mathcal{B}(\{k_{\lamdos}\})^2  \|f\|^2   + 2\delta(B)^{-2}\sum_j |c_j|^2 \nonumber\\
&\le C_2(r, \| \mu_0 \|_\ast , \delta(B))\sum_j |c_j|^2 , \ \label{tagdd}
\end{align}
where the last inequality holds by \eqref{tupi} and  because by Corollary \ref{muss}
the Bessel constant $\mathcal{B}\big(\{  k_{\lamdos} \}\big)^2$
has a bound that depends only on $r$ and $\| \mu_0 \|_\ast$.

\noi
So, if $\rho(\lambda_j, \lambda'_j) \le \eta  \le r$, \eqref{atla} and  \eqref{tagdd} yield
\begin{align}
\sum_j \Big| \sum_{i: i\ne j} \frac{c_i}{B_i(\lambda_i)}\,
\la B_{i} k_{\lambda_i},  k_{\lambda'_j} \ra \Big|^2
&=
\sum_j |\phi_{\lambda_j}(\lambda'_j)|^2  \, \Big| \sum_{i: i\ne j} \frac{c_i}{B_i (\lambda_i)}\,
\la B_{i,j} k_{\lambda_i},  k_{\lambda'_j} \ra \Big|^2   \nonumber \\
&=
\sum_j |\phi_{\lambda_j}(\lambda'_j)|^2    (1-|\lambda'_j|^2) | F_j(\lambda'_j)|^2 \nonumber\\
&\le \eta^2 \, C_3(r, \| \mu_0 \|_\ast , \delta(B))\sum_j |c_j|^2 .  \ \ \label{qu3}
\end{align}
Inserting inequalities \eqref{qu2}   and \eqref{qu3} in \eqref{qu1}, and using Cauchy-Schwarz again,
we obtain
\begin{align*}
\sum_j |\la f , k_{\lambda_j}-k_{\lambda'_j} \ra|^2 &\le
\eta^2 \,  2\!\left[\frac{4^2}{\delta(B)^2} +  C_3(r, \| \mu_0 \|_\ast , \delta(B))\right]
\sum_j |c_j|^2   .
\end{align*}
Since $c_j=\la f,k_{\lambda_j}\ra$, this proves the theorem for $f\in K_B$. A general $h\in H^2$ decomposes as $h= f +Bg$, where
$f\in K_B$, $g\in H^2$ and $\|h\|^2=\|f\|^2 + \|g\|^2$. Thus,
$$
\la f+Bg, k_{\lambda_j}-k_{\lambda'_j} \ra
=\la f, k_{\lambda_j}-k_{\lambda'_j} \ra - B(\lambda'_j) \la g, k_{\lambda'_j} \ra ,
$$
and since $|B(\lambda'_j)| = \rho(B(\lambda_j), B(\lambda'_j)) \le \rho(\lambda_j, \lambda'_j) \le \eta$,
\begin{align*}
\sum_j |\la f+Bg, k_{\lambda_j}-k_{\lambda'_j} \ra|^2 &\le
2\sum_j |\la f, k_{\lambda_j}-k_{\lambda'_j} \ra|^2 + 2\eta^2 \sum_j | \la g, k_{\lambda'_j} \ra|^2 \\
&\le
\eta^2 C(r, \| \mu_0 \|_\ast , \delta(B)) \, (\|f\|^2 + \|g\|^2),
\end{align*}
where  the last inequality uses the result for $f\in K_B$ and Corollary \ref{muss}.
\edem

\section{Separation conditions}
\noi
{\bf Definition.}
We say that a sequence (finite or not) $S$ in $\disc$ is $m$-separated (with radius $\ge \beta$) if
every pseudo-hyperbolic ball $\bru(z, \beta)$, with $z\in S$, has no more than $m$ points of
$S$ including repetitions.\\

\noi It is clear that if  we take $0<\beta_1 < \beta$ in the above definition then $S$ is also
$m$-separated with radius $\ge\beta_1$. Also, $(m-1)$-separated implies $m$-separated, and $1$-separated simply means
 separated, as in (2) of Lemma \ref{becar}.

Since the order of a sequence will not be relevant in what follows,
we operate with them as if they were sets with pointwise multiplicities.
So, for instance, the union of two sequences has the points of both with the sum
of multiplicities and some order.

\begin{lemma}\label{vx1}
Let $S$ be an $m$-separated sequence in $\disc$ with radius $\ge \beta$.
Then $S$ splits into at most $m$ separated sequences (finite or not) with radius $\ge \beta/4m$.
\end{lemma}
\bdem
Consider the balls $\bru(z_n, \beta/4m)$ including repetitions of the $z_n\in S$.
Let $U\subset \bigcup_{z_n\in S}  \bru(z_n, \beta/4m)$ be a connected component. Hence, $U$ is a union
of these balls  and we show that there cannot be more than $m$ of them (including repetitions). Otherwise,
$$U \supset \bru(z_{n(1)}, \beta/4m) \cup \ldots \cup \bru(z_{n(m+1)}, \beta/4m) ,
$$
where the union of balls is connected and has pseudo-hyperbolic diameter $\le 2(m+1) \beta/4m\le \beta$.
Consequently, $\bru(z_{n(1)}, \beta)$ contains the points
$z_{n(j)}$ for $j=1, \ldots , m+1$, contradicting the hypothesis.  Therefore
$$U = \bru(z_{n(1)}, \beta/4m) \cup \ldots \cup \bru(z_{n(k)}, \beta/4m) ,
\ \mbox{ where $k\le m$}.
$$
If we accept the empty set and finite sequences then $S = \bigcup_{j=1}^m S_j$, where
 $S_j = \{z_{n(j)}\in U: \, U \mbox{ is a connected componente}\}$  for $1\le j\le m$.
\edem

\begin{lemma}\label{vx2}
Let $S$ be an $m$-separated sequence in $\disc$.
Then there is a set $J\subset \{ 1, \ldots , m\}$ and parameters\/ $0< \eta_p < \gamma_p<1$, for $p\in J$,
such that $S$ splits into subsequences $S_p\, (p\in J)$ (finite or infinite), so that when $p\in J$:
\begin{enumerate}
\item[(i)] There is a single multiplicity sequence $S'_p \subset S_p$
such that
$$S_p = \bigcup \{ \bru( z'_n(p), \eta_p )\cap S : \ z'_n(p)\in S'_p\},$$
where each $\bru( z'_n(p), \eta_p )$ has $p$ points of $S_p$ counting
multiplicities.
\item[(ii)] $\rho (z'_n(p), z'_k(p)) > \gamma_p\ $ if $\ n\ne k$.
\item[(iii)] $\rho (S_p, S_k) > (4/5)\gamma_p\ $ if $\ p> k$, and\/ $p, k\in J$.
\item[(iv)]
Once $\gamma_p$ is obtained by reverse induction we can choose $0<\eta_p < \gamma_p$ arbitrarily small,
eventually lowering the index $p$ until $S_p \neq \emptyset$.
\end{enumerate}
\end{lemma}
\bdem
By hypothesis there is a radius $\beta_m >0$ such that $S$ is $m$-separated of radius $\ge \beta_m$.
Hence, by the previous lemma, $S$  splits into at most $m$ separated
sequences $T_1, \ldots, T_m$ of  separation $\ge \beta_m/4m: = \gamma_m$
(some could be empty or finite).
Chose any $0<\eta_m <  \gamma_m /10$,  and set
$$
S'_m = \{ z'_n \in T_m : \bru(z'_n , \eta_m) \mbox{ has just $m$ points of $S$ counting repetitions}\},
$$
and let $S_m$ be the sequence of all the points of $S$ in those balls.
If $S_m= \emptyset$, we reindex the parameters $\gamma_m$ and $\eta_m$ as
$\gamma_{m-1}$ and $\eta_{m-1}$, respectively, declare that $m\not\in J$ and
keep the process with $m-1$ instead of $m$.
Notice that by definition, $S_m=\emptyset\,$ if and only if $\,S'_m=\emptyset$, which holds if $\,T_m=\emptyset$.

If $S_m\neq \emptyset$, we keep $m\in J$ and notice that each of the balls has one and only one point
of each $T_q\ (1\le q\le m)$, because the distance between two different points in $T_q$ is $\ge \gamma_m$.
Thus,
\bequ\label{dory}
\rho (z'_n, z'_k) > \gamma_m \ \mbox{ if $n\ne k$}, \, \mbox{ and }\ \,
\rho(S_m, S\setminus S_m)  > \gamma_m - 2 \eta_m\ge \frac{4}{5}\gamma_m.
\eequ
If the remaining $S\setminus S_m$ is empty we are done. Otherwise it is a $(m-1)$-separated sequence
(with radius $\ge \eta_m/2$).
Indeed, if there is $z\in S\setminus S_m$ such that $\bru (z, \eta_m/2)$ has at least
$m$ points of $S\setminus S_m$ counting multiplicities, then this ball has at least one point
$z'$ of $T_m$, and consequently $\bru (z', \eta_m)$ has at least $m$ points
of $S$. In addition, since $\eta_m < \gamma_m/10<\beta_m$, it cannot have more than $m$ points of $S$.
Thus, $\bru (z', \eta_m)$ has exactly $m$ points of $S$, implying that
all the points of $S$ in $\bru (z', \eta_m)$ are in $S_m$, a contradiction.

Therefore, we can repeat the process above with $S\setminus S_m$ instead of $S$,
$m-1$ instead of $m$ and $\beta_{m-1}:= \eta_m/2$.
So, again there is $\gamma_{m-1}>0$ analogously defined and we can choose
$0<\eta_{m-1} <  \gamma_{m-1} /10 $, otherwise arbitrary, to define
analogous $S'_{m-1}$ and $S_{m-1}$ as before, just observing that
they could be empty, in which case $m-1\not\in J$ and we lower the index from $m-1$ to $m-2$.

We keep this process  going until we exhaust all the points of $S$.
Since the construction repeats condition \eqref{dory} for each $p\in J$, we get that if $p\in J$,
$$\rho (z'_n, z'_k) > \gamma_p \ \mbox{ for $z'_n, z'_k \in S'_p\ $ with $n\ne k$}
$$
and
$$
\rho(S_p, S\setminus {\textstyle \bigcup} \{ S_q: q\in J, \, q\ge p\} )  > \gamma_p - 2 \eta_p\ge \frac{4}{5}\gamma_p.
$$
Thus, if $p, k\in J$, with $p>k$, then $S_k\subset S\setminus \bigcup \{ S_q: q\in J, \, q\ge p\}$,
and consequently
$$
\rho(S_p, S_k ) \ge \frac{4}{5}\gamma_p.
$$
Therefore (ii) and (iii) hold. Also,  (i) and (iv) hold by construction.
\edem

\vspace{1mm}\noi
Suppose that $S=\{ \lambda_j\}$ is $m$-separated and $\sum (1-|\lambda_j|^2) \delta_{\lambda_j}$ is a
Carleson measure. By Lemma \ref{vx1} and  Lemma \ref{becar}, $S$ is the union of at most
$m$ interpolating sequences. More importantly for our purpose,
each $S_p = \{ \lambda_j(p) : \, j\ge 1 \}$ of the decomposition
given by Lemma \ref{vx2} is a finite union of interpolating sequences.
Let $B_p$ be the Blaschke product whose zeros are $S_p$ counting multiplicities,
or $B_p\equiv 1$ if $S_p=\emptyset$.
It is then known (see \cite{kerr}) that
\bequ\label{vchi}
\inf\{ |B(z)| : \rho(z, S_p) \ge \gamma \}>0 \, \ \mbox{ for any } \gamma>0.
\eequ
That is, $B_p$ is bounded below away from zero at any fixed positive distance of its zeros,
which fails for arbitrary Blaschke products.  Now define
$$
A_\nuq := \prod_{p=1, \, p\neq \nuq}^m  B_p
\ \ \mbox{ and }\ \
f_\nuq^i = \sum_j c_j(\nuq) \alpha_j^i(\nuq) k_{\lambda_j(\nuq)} \in K_{B_\nuq},
$$
for $1\le i, \nuq\le m$\/ and $\{c_j(\nuq)\}_j\in \ell^2$,
where \/$|\alpha_j^1(\nuq)|^2 +\cdots + |\alpha_j^m(\nuq)|^2=1$
for all $1\le \nuq \le m\,$ and  $j\ge 1$. Also, notice that $f^i_\nuq=0= c_j(\nuq)$ if $B_\nuq \equiv 1$.

Since the zeros of $A_\nuq$ are $\bigcup_{p=1, \, p\neq \nuq}^m  S_p $ and by (iii) of Lemma \ref{vx2},
$$\rho \Big(\bigcup_{p=1, \, p\neq \nuq}^m  S_p , S_\nuq\Big)>0,
$$
it follows from \eqref{vchi}  that
$$0 < \varepsilon_\nuq = \inf \{ |A_\nuq(\lambda_j(\nuq))|:   \, j\ge 1  \}.
$$
We also need an elementary fact about Toeplitz operators.
If $g\in \papa(\disc)$, it is easy to prove that
the normalized reproducing kernels for $H^2$ are eigenvalues of the Toeplitz operator
$T_{\ov{g}}$ such  that $T_{\ov{g}} k_\lambda = \ov{g(\lambda)} k_\lambda$. \\

\noi
The next proposition will allow us to reduce the proof of sufficiency of Theorem \ref{kasem} to
a particular  case. We keep the above notations.

\begin{propos}\label{blakees-m}
If
\bequ\label{lox}
\sum_{i=1}^m \| f_\nuq^i \|^2 \ge  D_\nuq ^2 \sum_j |c_j(\nuq)|^2
\eequ
for all\/ $1\le \nuq\le m$ and $\{c_j(\nuq)\}_j \in \ell^2$, then
$$
m \sum_{i=1}^m \| f_1^i + \cdots + f_m^i \|^2 \ge
\min_{1\le \nuq\le m} \{D_\nuq^2 \varepsilon_\nuq^2\} \sum_{\nuq=1}^m    \sum_j |c_j(\nuq)|^2 .
$$
\end{propos}
\bdem
Since $T_{\ov{A}_\nuq}$ is a contraction on $H^2$, for each $1\le i\le m$,
 \begin{align*}
 \| f_1^i + \cdots + f_m^i \|^2  &\ge \|  T_{\ov{A}_\nuq} (f_1^i + \cdots + f_m^i) \|^2
 \, =\,  \|  T_{\ov{A}_\nuq} f_\nuq^i \|^2  \\*[1mm]
 &= \big\| \sum_j c_j(\nuq) \alpha_j^i(\nuq) \ov{A_\nuq(\lambda_j(\nuq))} k_{\lambda_j(\nuq)} \big\|^2
 \end{align*}
for $1\le \nuq\le m$. So, by \eqref{lox},
$$
\sum_{i=1}^m  \| f_1^i + \cdots + f_m^i \|^2 \ge D_\nuq^2 \sum_j |c_j(\nuq)|^2 \, |A_\nuq(\lambda_j(\nuq))|^2
 \ge D_\nuq^2 \sum_j |c_j(\nuq)|^2 \varepsilon_\nuq^2 .
$$
Adding these inequalities for $1\le \nuq\le m$, we obtain
$$
 m\sum_{i=1}^m \| f_1^i + \cdots + f_m^i \|^2   \ge \sum_{\nuq=1}^m  D_\nuq^2 \varepsilon_\nuq^2 \sum_j |c_j(\nuq)|^2
  \ge   \min_{1\le \nuq\le m} \{D_\nuq^2 \varepsilon_\nuq^2\} \, \sum_{\nuq=1}^m \sum_j |c_j(\nuq)|^2 .
\vspace{-5mm}
$$
\edem

\subsection{The sufficiency of Theorem \ref{kasem}}

\noi Finally, the proof of Theorem \ref{kasem} needs the following elementary inequality.

\begin{lemma}
If $x, y_1 , \ldots , y_p \in H^2$, with $p\le m$, then
\bequ
\Big\| x+  \sum_{k=1}^p y_k \Big\|^2  \ge \frac{\| x \|^2}{2}  - m\sum_{k=1}^p      \| y_k \|^2 .
\label{ichito}
\eequ
\end{lemma}
\bdem
This follows  using Cauchy-Schwarz's inequality twice. For $y\in H^2$,
$\| x \|^2   \le 2 \| x+y \|^2 +  2 \| y \|^2$,
so $\| x+y \|^2 \ge \| x \|^2/2  - \| y \|^2$.
Then,
$$
\Big\| x+  \sum_{k=1}^p y_k \Big\|^2 \ge \frac{\| x \|^2}{2}  - \Big\| \sum_{k=1}^p y_k \Big\|^2
 \ge \frac{\| x \|^2}{2}  - p\sum_{k=1}^p      \| y_k \|^2 .
\vspace{-4mm}
$$
\edem

\bdem[Proof of sufficiency for Theorem \ref{kasem}]
\hspace{1mm} Since by (1) of the theorem, $S=\{\lambda_j\}$ is $m$-separated and
$\mu= \sum_{i}  (1-|\lambda_i|^2) \delta_{\lambda_i}$ is a Carleson measure,
Lemma \ref{vx2} and the comments that follow the lemma apply to $S$.
Therefore, we have the decomposition  of the lemma
$$
S= \bigcup_{p\in J} S_p, \ \mbox{ with }  J\subset \{1,\ldots , m\},
$$
where by (iv) we can choose $\eta_p<\eta$ for all $p\in J$
(here $\eta$ is the parameter of Theorem \ref{kasem}).

This  guarantees that the sequences $S_p$ satisfy (2) of the theorem for each
ball appearing in (i) of Lemma \ref{vx2}.

Furthermore, Proposition \ref{blakees-m} reduces the problem of proving \eqref{mini} for the sequence $S$
to prove it for each sequence $S_p$ from the above decomposition.
Therefore we fix an arbitrary $p$ between $1$ and $m$.
The subsequence $S'_p$ in (i) of Lemma \ref{vx2} is
separated with radius of separation $\ge \gamma_p$, so Lemma \ref{becar} says that
it is itself interpolating.

We re-index $S_p$ as follows: write $S'_p= \{ \lambda_j(1) : j\ge 1 \}$ and
$\{ \lambda_j(\nu) : 1\le \nu\le p\}$ for the $p$ elements of $S_p$ that are in
$\bru( \lambda_j(1) , \eta_p)$. Therefore
$$
S_p = \{ \lambda_j(\nu) : 1\le \nu\le p, \ 1\le j\} ,
$$
and the respective parameters of Theorem \ref{kasem} write as
$c_j(\nu)$ and $\alpha_j^i(\nu)$, for $1\le \nu\le p, \ 1\le j$ and $1\le i\le m$.
We see that for each $1\le i\le m$,
\begin{eqnarray}
\lefteqn{  \big\|  \sum_{j\ge 1} \sum_{\nu=1}^p   c_j(\nu) \alpha_j^i(\nu) k_{\lambda_j(\nu) }\big\|^2 = }
\nonumber \\
&=& \hspace{-1mm}
\Big\|  \sum_{j\ge 1} \Big[
\Big(\sum_{\nu=1}^p  c_j(\nu) \alpha_j^i(\nu) \Big)\,  k_{\lambda_j(1)} +
\sum_{\nu=2}^p  c_j(\nu) \alpha_j^i(\nu)\, \big(k_{\lambda_j(\nu)}- k_{\lambda_j(1)}\big)\Big] \Big\|^2
\nonumber \\
&\hspace{-4mm}\mbox{}&
\hspace{6mm}
\ge \, \frac{1}{2} \Big\|
\sum_{j\ge 1} \Big(\sum_{\nu=1}^p  c_j(\nu) \alpha_j^i(\nu) \Big)\, k_{\lambda_j(1)} \Big\|^2 -
\nonumber \\
&\ &
\hspace{6mm}
-\  m \sum_{\nu=2}^p  \Big\|  \sum_{j\ge 1}  c_j(\nu) \alpha_j^i(\nu)\,
\big( k_{\lambda_j(\nu)}- k_{\lambda_j(1)}\big)  \Big\|^2,
\label{jd0}
\end{eqnarray}
where the inequality holds by \eqref{ichito}.
If $B$ is the Blaschke product with zeros $\{\lambda_j(1)\}_j$, Lemma \ref{becar} says that
$\delta(B)$ is estimated depending on $\gamma_p$  and the measure $\mu$ associated to $S$.
Thus, the comments that follow  \eqref{ddta}
imply that
$\{ k_{\lambda_j(1)}: \, j\ge 1\}$ is a Riesz sequence with lower constant
$D_1^2$ independent of $\eta_p$ (only depends on $\gamma_p$  and $\mu$).
Hence,
$$
\big\|
\sum_{j\ge 1} \big(\sum_{\nu=1}^p  c_j(\nu) \alpha_j^i(\nu) \big)\,  k_{\lambda_j(1)}  \big\|^2
\ge D_1^2
\sum_{j\ge 1} \big|\sum_{\nu=1}^p  c_j(\nu) \alpha_j^i(\nu) \big|^2,
$$
and taking into account that, with the new indexation, the last condition of the theorem rewrites as
$$
D_0^2
\left\|
    \begin{bmatrix}
      c_j(1) \\
    \vdots \\
    c_j(p) \\
  \end{bmatrix}
\right\|^2_{\mathbb{C}^p}
\le
\left\|
  \begin{bmatrix}
    \alpha_j^1(1)  & ...  &\alpha_j^1(p) \\
    \vdots     &            &  \vdots  \\
     \alpha_j^m(1)  & ...  &\alpha_j^m(p)\\
  \end{bmatrix}
  \begin{bmatrix}
    c_j(1) \\
    \vdots \\
    c_j(p) \\
  \end{bmatrix}
\right\|^2_{\mathbb{C}^{m \times 1}}\hspace{-8mm},
\ \ \ \ \ \ \ \ \forall
 \begin{bmatrix}
    c_j(1) \\
    \vdots \\
    c_j(p) \\
  \end{bmatrix}
\in \mathbb{C}^p,
$$
we obtain
\bequ
\sum_{i=1}^m \big\|
\sum_{j\ge 1} \big(\sum_{\nu=1}^p  c_j(\nu) \alpha_j^i(\nu) \big)\,  k_{\lambda_j(1)}  \big\|^2
\ge D_1^2 D_0^2
\sum_{j\ge 1} \sum_{\nu=1}^p   \big|c_j(\nu)  \big|^2. \label{jd1}
\eequ
Observe that, as with $D_1$, the constant $D_0$ does not depend on $\eta_p$,
since it only depends on the $\eta\,$ in condition (2) of Theorem \ref{kasem}.

To estimate the second term of \eqref{jd0} we notice that since $\eta_p<\eta$, by  Theorem \ref{difere}
there is a constant $C_S>0$ depending only on the original sequence $S$ (and therefore is independent of $\eta_p$), such that
\begin{align}
\big\|  \sum_{j\ge 1}  c_j(\nu) \alpha_j^i(\nu)\, \big( k_{\lambda_j(\nu)}- k_{\lambda_j(1)}\big)\big\|^2
&\le
\mathcal{B}\big( \{k_{\lambda_j(\nu)}- k_{\lambda_j(1)}\}\big)^2\sum_{j\ge 1}  \big| c_j(\nu) \alpha_j^i(\nu)\big|^2
\nonumber \\
&\le C_S^2 \eta_p^2\sum_{j\ge 1}  \big| c_j(\nu) \alpha_j^i(\nu)\big|^2  ,\nonumber
\end{align}
for all $2\le \nu \le p$ and $1\le i\le m$. Thus
\bequ
\sum_{i=1}^m \sum_{\nu=2}^p \big\|  \sum_{j\ge 1}
c_j(\nu) \alpha_j^i(\nu)\, \big( k_{\lambda_j(\nu)}- k_{\lambda_j(1)}\big)\big\|^2
\le C_S^2 \eta_p^2 \sum_{\nu=2}^p \sum_{j\ge 1}  \big| c_j(\nu)\big|^2  .\label{jd2}
\eequ
If we insert \eqref{jd1} and \eqref{jd2} into \eqref{jd0}, we get
\begin{eqnarray*}
\lefteqn{ \sum_{i=1}^m \Big\|  \sum_{j\ge 1} \sum_{\nu=1}^p
c_j(\nu) \alpha_j^i(\nu) k_{\lambda_j(\nu) }\Big\|^2 \ge }   \nonumber \\
&\ge&
\frac{1}{2} D_1^2 D_0^2
\sum_{j\ge 1} \sum_{\nu=1}^p   \big|c_j(\nu)  \big|^2 -
m   C^2_S \eta_p^2 \sum_{\nu=2}^p \sum_{j\ge 1}  \big| c_j(\nu)\big|^2 \\
&\ge&
\Big(\frac{1}{2} D_1^2 D_0^2 - m   C^2_S \eta_p^2 \Big)
\sum_{j\ge 1} \sum_{\nu=1}^p   \big|c_j(\nu)  \big|^2 .
\end{eqnarray*}
We can take $\eta_p$ small enough for the expression between brackets to be $\ge (1/4)  D_1^2 D_0^2$,
which proves the theorem.
\edem

\noi
We go back to our original problem.
By \cite[Thm.\@ 5.6]{accmp} any normal operator $N$ on a Hilbert space that admits
a finite union of orbits as a frame is diagonalizable.
Therefore $N$ is unitarily equivalent to an operator $A$ as in:

\begin{theo}
Let $A$ be a diagonal operator with respect to the standard basis
in $\ell^2(J)$, where $J= \mathbb{N}$ or it is finite, and let $a^1, \ldots , a^m \in \ell^2(J)$.
Then $\{ A^n a^i: \, n\in \mathbb{N}_0, \, 1\le i\le m\}$ is a frame if and only if
\begin{itemize}
\item Each $a^i$   is given by \eqref{niunii} for $\,i=1, \ldots , m$.
\item The sequence of eigenvalues $\{ \lambda_j : \, j\in J\}$ of $A$ and the double sequence
$\{ \alpha^i_j : \, j\in J, \,1\le i\le m\}$ appearing in \eqref{niunii} satisfy Theorem \ref{kasem}.
\end{itemize}
\end{theo}

\noi
We finish the paper with a comment on the tails of the orbits and further remarks,
keeping the previous notations and assumptions. If $n_0\ge 1$ is an integer,
the calculation leading to \eqref{kcat} gives

$$
\sum_{n\ge n_0} |\la A^n \tilde{a}^i , c \ra|^2 =
\big\|   \sum_{j}\alpha^i_j c_j \ov{\lambda}^{n_0}_j   k_{\lambda_j}\big\|^2
\ \ \ \mbox{ for $1\le i\le m$}.
$$
Therefore, if \eqref{mini} holds and we write $\mathcal{I}=\{ j: \lambda_j \neq 0  \}$,
$$
 \sum_{i=1}^m   \sum_{n\ge n_0} |\la A^n \tilde{a}^i , c \ra|^2 \ge
 D^2 \, \sum_{j}   |\ov{\lambda}^{n_0}_j  c_j|^2 \ge
 D^2\, \min_{j\in \mathcal{I}} \{|\lambda_j|^{n_0}\}  \, \sum_{j\in \mathcal{I}}   |c_j|^2,
$$
for all $c\in \ell^2$, where $\lambda_j \not \to 0$ by  (1) of Theorem \ref{kasem}.
This means that $\{ A^n \tilde{a}^i : \, n\ge n_0, \, 1\le i\le m\}$ is a frame
for $(\mbox{Ker}\, A)^\bot$.\\

\noi
In \cite[Thm.$\,$3.4]{CHP} the authors characterize the bounded operators $T$ on a separable Hilbert space
of infinite dimension $\mathcal{H}$
(and the vectors $f\in \mathcal{H}$) such that $\{ T^n f: \, n\in\N\cup\{0\}\}$ is a frame, as
those that are similar to the forward shift on $H^2$ or its compression to
one of the spaces $K_u$, where $u$ is an inner function other than a finite Blaschke product.
As said before, a more precise statement was given in \cite{acmt} when $T$ is diagonal,
and  the case of a normal operator $T$ was reduced to the diagonal case in \cite{accmp}.
This result is also deduced in \cite{CHP} as a consequence of their general theorem and known facts about model spaces.
Techniques from the theory of model spaces are also used for studying frames and Bessel sequences in
\cite{H-J-W} and \cite{Jorge-Ma}.
Theorem \ref{kasem} here is more related to a particular case of model spaces $K_u$, in which $u$ is a finite
product of interpolating Blaschke products.

\end{document}